\documentclass[british,a4paper]{amsart}
\usepackage{babel}
\usepackage[T1]{fontenc}
\usepackage{amscd}
\usepackage{url}
\usepackage{amssymb}

\newtheorem{theorem}{Theorem}[section]
\newtheorem{corollary}[theorem]{Corollary}
\newtheorem{proposition}[theorem]{Proposition}
\newtheorem{lemma}[theorem]{Lemma}

\theoremstyle{remark}
\newtheorem*{remark}{Remark}

\theoremstyle{definition}
\newtheorem{definition}{Definition}[section]

\newtheorem{example}{Example}

\newcommand{\ie}{\textit{i.e. }}

\newcommand{\p}{\mathbb{P}}
\newcommand{\F}{\mathcal{F}}

\newcommand{\T}{\mathbb{T}}

\newcommand{\Gr}{\mathbb{G}}

\newcommand{\OO}{\mathcal{O}}
\newcommand{\ud}{\mathrm{d}}

\newcommand{\abs}[1]{\lvert#1\rvert}

\newcommand{\Hh}{\mathcal{H}}

\DeclareMathOperator{\rk}{rk}

\DeclareMathOperator{\Hilb}{Hilb}

\DeclareMathOperator{\length}{length}
\DeclareMathOperator{\Al}{Al}
\DeclareMathOperator{\sing}{sing}
\DeclareMathOperator{\red}{red}

\DeclareMathOperator{\Hom}{Hom}

\begin{document}


\title[Congruences in $\p^4$ with irreducible fundamental Surface]{On first 
order Congruences of Lines in $\p^4$ with irreducible fundamental Surface}

\author{Pietro De Poi}
\thanks{This research was partially supported by the DFG Forschungsschwerpunkt
``Globalen Methoden in der Komplexen Geometrie'', and the EU (under the EAGER
network).}
\address{Mathematisches Institut\\
Universit\"at Bayreuth\\
Lehrstuhl VIII \\
Universit\"atsstra\ss e 30\\
D-95447 Bayreuth\\
Germany\\
}

\begin{abstract}
In this article we study congruences of lines in $\p^n$, 
and in particular of order one. After giving general results,  
we obtain a complete classification in the case of $\p^4$ 
in which the fundamental surface $F$ is in fact a variety---\ie it is 
integral---and the congruence is the irreducible set of the trisecant lines 
of $F$. 
\end{abstract}

\maketitle                   

\section{Introduction}
A congruence of lines in $\p^n$ is a family of lines of 
dimension $n-1$, and its order is the number of lines passing 
through a general point in $\p^n$. Here we are interested in classifying 
congruences of order one in $\p^4$ with reduced and 
irreducible \emph{fundamental surfaces}.  In this situation, this means  
that the congruence is the closure of the trisecant line of its fundamental surface, 
see below. 
 
The first thing one can say, about first order congruences, is that---thanks to the 
Zariski main theorem---the focal 
locus $\Phi$ coincides with the fundamental locus, \ie it is characterized by the 
fact that through a point in it there pass infinitely many lines of the 
family. Although for a general congruence in $\p^n$ one expects that 
$\Phi$ is---since it is the image of the ramification divisor---of dimension $n-1$, 
in the case in which we are interested, 
since it coincides with the fundamental locus, usually it is reducible 
(and non-reduced) of dimension (at most) $n-2$. More precisely, the congruence is 
(a subvariety of) the family of the $(n-1)$-secant lines to the focal locus $\Phi$, 
see Proposition~\ref{prop:fofi} below. 
In particular, in the case of first order congruences in $\p^4$, $\Phi$ 
is, in general, a reducible and non-reduced scheme of dimension two; 
see~\cite{DP2}, Example~\textup{(2)}. 
Here we are interested in the case in which the pure focal locus is in fact a 
variety (\ie it is reduced and irreducible) of dimension two, 
and the congruence
is given by the family of the trisecant lines of this fundamental 
surface---provided that this family is irreducible. 
For the history, other general results and references about first order 
congruences, see~\cite{PDP}, \cite{DP2} and \cite{DP3}.

In this paper the study of these congruences is carried on: 
we extend the classification of 
Ascione, Severi and Marletta (see~\cite{A}, \cite{S} and \cite{M1}, see also 
\cite{Au}), of the first order congruences in $\p^4$ given by the family of 
the trisecant lines of a surface $F$---provided that the family is 
irreducible---without assumptions on the singularities of $F$. 
It is interesting that---in this general case---the order of the family of the 
trisecant lines of an irreducible surface $F$ in $\p^4$ is less than or 
equal to the number of its apparent triple points, while in the classical 
studies it was always equal, since they allowed only general singularities 
for the surface $F$ (see for example in \cite{M1}). We will also show, in 
Example~\ref{ex:1}, a surface without trisecants but with one apparent triple 
point. We note here that some of the surfaces, or, more generally, varieties 
with one apparent triple point can be constructed---as observed in Remark 2 
of \cite{AR}---via the so called completely \emph{subhomoloidal systems} of hypersurfaces 
(see Definition 1 of \cite{AR}); actually, in this case of cubic hyperplanes. 
In fact, their construction is a generalization of the one given---obtained 
independently from that of \cite{AR}---in (the proof of) Theorem~2.6 of \cite{DP3}. 
Concerning the trisecant lines to a surface in $\p^5$, there are results of 
I. Bauer, in \cite{Icr} and \cite{Ito} (see also \cite{FR}).

This article is structured as follows: 
after giving, in Section~\ref{sec:1}, the 
basic definitions, we consider some general results about 
congruences in $\p^n$, which are in fact an extension of the results of 
\cite{DP3}. We start with 
the calculus of the degree of the hypersurface in $\p^n$ which corresponds to  
a (special) linear section of a congruence. 
We use this in Subsection~\ref{sec:3} to study the components of the 
intersection of two of these hypersurfaces.
After giving the central definition of fundamental $d$-loci for congruences 
of order one, we
study in Subsection~\ref{sec:par} 
all the other subschemes of the focal locus of this situation, which we will 
call the parasitic schemes. 
This is a generalization of the classical definition of the parasitic 
planes for a first order congruence in $\p^4$ (see \cite{A} and 
\cite{M1})---and in this last case of first order congruences in $\p^4$ 
the two definitions agree. 
In this subsection the first 
important properties of these schemes are proved; it is worth noting that 
the family of lines contained in a parasitic scheme of dimension $d$ 
has dimension at least $d$.
We finish this section by giving some general formulae from which we get a 
degree bound for the fundamental $(n-2)$-locus---if it coincides with the 
pure fundamental locus. 
Unfortunately, this bound is not sharp, as we will see in the case in $\p^4$. 

The congruences of order one in $\p^4$ with an irreducible fundamental surface 
are classified in 
Section~\ref{sec:1t}. 
We obtain, in Theorem~\ref{thm:primo}, a complete list (for further details about these 
congruences, the fundamental 
surfaces and the parasitic planes, see Section~\ref{sec:1t}); for stating this theorem, 
we need to say that for us a---possibly singular---\emph{Del Pezzo surface in $\p^5$}  
is as defined in \cite{dem}, \ie $\p^2$ blown-up in $4$ points in 
\emph{quasi-general position} (``position presque generale'' of \cite{dem}, pag. 39), embedded via its 
anticanonical system; in particular the
surface is normal, see pag. 63 of \cite{dem}, and so it has at most isolated singularities,  
while a---possibly singular---\emph{Bordiga surface in $\p^4$} is the degeneracy locus of a map
$\phi\in\Hom(\OO_{\p^4}^{\oplus 3},\OO_{\p^4}^{\oplus 4}(1))$---if 
the minors vanish in the expected (irreducible) codimension two. We will see in 
Lemma~\ref{lem:bor} that such a surface has at most five singular points. 
A classical reference for the Bordiga surfaces is in \cite{room}, Chapter~XIV.

\begin{theorem}\label{thm:primo} 
The irreducible surfaces in $\p^4$ whose trisecant lines generate an 
irreducible first order congruence are 
\begin{enumerate}
\item\label{ve} a smooth (projected) Veronese surface; 
\item\label{dp} a projection of a---possibly singular---Del Pezzo surface
in $\p^5$ from a point not lying on a plane containing an irreducible conic 
of the surface; 
\item\label{sc} a projection of a smooth rational normal scroll $S$ 
in $\p^6$ from a line 
which does not intersect $S$;  
\item\label{bo} a Bordiga surface (which has at most five singular points). 
\end{enumerate}

\end{theorem}

We finish the article studying which of the congruences of 
Theorem~\ref{thm:primo} 
are smooth, if interpreted as $3$-folds in the Grassmannian $\Gr(1,4)$. 
For simplicity, we restrict ourself to the cases of smooth surfaces.

\begin{theorem}\label{thm:smo}
The smooth surfaces in $\p^4$, whose trisecant 
lines form an irreducible first order congruence $B$ 
which is smooth in $\Gr(1,4)$, are
\begin{enumerate}
\item the smooth projected Veronese surface, in which case $B$ is a linear 
congruence, 
which has bidegree $(1,2)$ and sectional genus one; it is a 
Fano variety of index two (for more details, see 
Subsection~\textup{2.1} of \textup{\cite{DP3}}); 
\item the smooth Bordiga surface, in which case $B$ has bidegree 
$(1,8)$ and sectional genus $10$ (see Subsection~\textup{
2.2} of \textup{\cite{DP3}}
for more details).
\end{enumerate}
\end{theorem}
We note finally that linear congruences in $\p^4$ were studied 
classically---among others---by G. Castenuovo in \cite{En}, and G. Fano, 
in \cite{Fa} and \cite{Fa1} which are the source of inspiration for 
the celebrated works on Fano varieties by V. A. Iskovskih (\cite{I1} and \cite{I2}). 

\section{Notations, definitions and general results}\label{sec:1}
We will work with schemes and varieties over the complex field $\mathbb C$.
with the standard definitions and notation as in \cite{H}.
For us a \emph{variety} will always be irreducible and projective.
More information about general results and references 
about families of lines, focal diagrams and congruences can be found in 
\cite{DP2} or \cite{PDP}. Besides, we refer to \cite{GH} for notations about 
Schubert cycles and to \cite{Fu} for the definitions and results of 
intersection theory. So we denote the general Schubert cycle by: 
$$
\sigma_{a_0,a_1}:=[\{\ell\in\Gr(1,n)\mid a_0\ge a_1, \ell\cap\p^{n-1-a_0}\neq \emptyset, 
\ \ell,\p^{n-1-a_0}\subset\p^{n-a_1}\}]. 
$$
Here we recall that a \emph{congruence of lines} in $\p^n$ 
is a flat family $(\Lambda,B,p)$ of lines in $\p^n$ obtained by the 
desingularization of a subvariety $B'$ of dimension $n-1$ of the Grassmannian 
$\Gr(1,n)$ of lines in $\p^n$. 
We can summarize all the notations which are used in the following two diagrams: 
the first one defines the family   
$$
\begin{CD}
\Lambda:=\psi^*(\Hh_{1,n}) @>\psi^*>> 
\Hh_{1,n}@>p_2>>\p^n\\
 @VpVV  @Vp_1VV \\
 B@>\psi >>B'\subset\Gr(1,n),
\end{CD}
$$
where $\Hh_{1,n}\subset\Gr(1,n)\times \p^n$ is the incidence variety 
and $\psi$ is the 
desingularization map, and the second one explains the notation for the 
elements of the family 
$$
\begin{CD}
\Lambda_b\subset\Lambda@>f:=p_2\circ \psi^*>> 
\p^n\supset\Lambda(b):=f(\Lambda_b)\\
 @VpVV   \\
b\in B.
\end{CD}
$$

So $p$ is the restriction of the first projection 
$p_1:B\times\p^n\rightarrow B$ to $\Lambda$, while we will denote the 
restriction of the second projection by $f:\Lambda \rightarrow\p^n$. 
$\Lambda_b:=p^{-1}(b)$, $(b\in B)$ will be a line of the family and 
$f(\Lambda_b)=:\Lambda(b)$ is a line in $\p^n$. From the definition, we 
get that $\Lambda$ is a smooth variety of dimension $n$.

A point $y\in \p^n$ is called 
\emph{fundamental} if its fibre $f^{-1}(y)$ has dimension greater than the 
dimension of the general one. The \emph{fundamental locus} is the set of the 
fundamental points. It is denoted by $\Psi$.
The \emph{subscheme of the foci of the first order} or, simply, the 
\emph{focal scheme} 
$V\subset \Lambda$ is the scheme of the ramification points of $f$.
The \emph{locus of the first order foci}, or, simply, the \emph{focal locus}, 
$\Phi:=f(V)\subset \p^n$, is the set of the
branch points of $f$. In this article, as we did in \cite{DP2}, 
we will endow this locus with the 
scheme structure given by considering it as the scheme-theoretic image of 
$V$ under $f$ (see, for example, \cite{H}).

To a congruence is associated a \emph{sequence of degrees} 
$(a_0,\dotsc,a_\nu)$ if we write 
$$
[B]=\sum_{i=0}^{\nu}a_i\sigma_{n-1-i,i}
$$ 
---where we put $\nu:=\left [\frac{n-1}{2}\right ]$ and $[B]$ denotes the rational 
equivalence class of $B$---as 
a linear combination of 
Schubert cycles of the Grassmannian; in particular, the \emph{order} $a_0$ 
is the number of lines of $B$ passing through a general point in $\p^n$, and 
the \emph{class} $a_\nu$ is the number of lines intersecting a general 
$\nu$-plane and contained in a general $(n-\nu)$-plane 
(\ie as a Schubert cycle, $[B]\cdot \sigma_{n-1-\nu,\nu}$).
The fundamental locus is contained in the focal locus and the two loci 
coincide in the case of a first order congruence, see Proposition~\ref{prop:seg} \ie through a focal 
point there will pass infinitely many lines of the congruence.
An important result---independent of order and class---is the following:
\begin{proposition}\label{prop:fofi}\textup{ (C. Segre, \cite{sg}).}
On every line $\Lambda_b\subset \Lambda$ of the family, 
the focal subscheme $V$ 
either coincides with the whole  $\Lambda_b$---in which case $\Lambda(b)$ 
is called \emph{focal line}---or is a zero dimensional 
subscheme of $\Lambda_b$ of length $n-1$. Moreover, in the latter case, 
if $\Lambda$ is a first order congruence, $\Phi\cap \Lambda(b)$ has length  
$n-1$. 
\end{proposition}

See \cite{DP2} for a proof. 


\subsection{Linear sections of a congruence in $\Gr(1,n)$}\label{sec:lsc} 
Let $B(\subset\Gr(1,n))$ be a congruence with sequence of degrees 
$(a_0,\dotsc,a_\nu)$;
then we have that

\begin{proposition} \label{prop:ai}
Let $V_\Pi$ be the scroll
given by the lines of the congruence which 
meet a (fixed) general
$(n-2)$-plane $\Pi$. Then $V_\Pi$ is a hypersurface in $\p^n$ of degree 
$a_0+a_1$.  
\end{proposition}

\begin{proof}
We denote 
$p(f^{-1}(V_\Pi))=:G_\Pi$; clearly, as a Schubert cycle,  
$[G_\Pi]=[B]\cdot \sigma_1$, and an easy application of Pieri's formula gives
\begin{equation}\label{eq:bsi}
[G_\Pi] =\sum_{k=1}^{\nu}((\sum_{i=0}^k a_i)\sigma_{n-k,k}).
\end{equation}
Since $V_\Pi$ is a hypersurface, to obtain its degree we can intersect it
with a general line and compute the length of the zero dimensional subscheme
so obtained, or, which is the same, calculate the intersection of 
$[G_\Pi]$ with the Schubert cycle
$\sigma_{n-2}$, 
\ie the lines which meet a general line. 
So, from the intersection formula
for complementary Schubert cycles,
we have the claim. 
\end{proof}

\begin{remark}
If $B$ is a first order congruence, $\ell$ is a line of it
not contained in $V_\Pi$ and $P$ is a point in $V_\Pi\cap \ell$,
then $P$ is a focus for $B$, since at least two lines of the 
congruence pass through it.
\end{remark}


\subsection{Congruences of lines in $\p^n$ with $\Phi=\Psi$}
\label{sec:3}
From now on we will consider congruences such that the fundamental locus 
coincides with the focal locus (\ie $\Phi=\Psi$). We will see that with this 
hypothesis a congruence has order either zero or one (and \textit{vice versa}). 
First of all we see that  

\begin{proposition}\label{prop:seg}[C. Segre, \cite{sg}]
The fundamental locus of a congruence $\Lambda$ of order either zero or one 
coincides with the focal locus and has dimension at most $n-2$.
\end{proposition}

\begin{proof} 
The case of order $a_0=0$ is immediate by dimensional reasons. Let us now 
consider the case  $a_0=1$: 
the fact that the two loci coincide is a straightforward 
consequence of the Zariski Main Theorem, see \cite{H},  
since the map $f$ is generically $(1:1)$. 

The fundamental locus $F$ cannot have dimension $n-1$;
otherwise the subscheme of the first
order foci $V$ would coincide with $\Lambda$, and this would 
contradict the fact that we have a $(1:1)$ map.
\end{proof}

The following result was suggested to us by F. Catanese: 

\begin{theorem}\label{thm:cata}
Let $B$ be a congruence whose focal locus $F$ has codimension at least two.
Then, $B$ has order either zero or one.
\end{theorem}

This theorem can be found in \cite{DP1}. The proof is as follows:  
consider the restriction of the map $f:\Lambda\rightarrow\p^n$ to the
set $\Lambda\setminus f^{-1}(F)$. Then, either  $f^{-1}(F)=\Lambda$, in which
case $B$ is a congruence of order zero, or the map $f\mid_{\Lambda\setminus f^{-
1}(F)}$
defines an unramified covering of the set $\p^n\setminus F$. But it is a
well-known fact that---by dimensional reasons---$\p^n\setminus F$ is
simply connected and $\Lambda\setminus f^{-1}(F)$ is connected. Therefore,
$f\mid_{\Lambda\setminus f^{-1}(F)}$ is a homeomorphism, hence $f$ is a 
birational map and $B$ is a first order congruence.

\begin{corollary}
A congruence in $\p^n$ has order zero or one if and only if the focal locus $\Phi$ coincides with 
the fundamental locus $\Psi$; moreover, $\dim(\Phi)\le n-2$. 
\end{corollary}

\begin{proposition}\label{thm:pp}
Let $\Lambda$ be a congruence such that $\Phi=\Psi$, and let $\Pi$ and 
$\Pi'$ be general $(n-2)$-planes in $\p^n$. Then the complete intersection
of the hypersurfaces  $V_\Pi$ and  $V_{\Pi'}$ is a (reducible) 
$(n-2)$-dimensional scheme  $\Gamma$ which---set-theoretically---is the 
union of 
the $(n-2)$-dimensional scroll $\Sigma$ given by the lines of the congruence 
meeting $\Pi$ and $\Pi'$, which has degree $a_0+2a_1+a_2$ if $n>3$ and 
$a_0+a_1$ for $n\le 3$ and the focal locus $\Phi$ (possibly as an embedded 
component of $\Sigma$).
\end{proposition}

\begin{proof}
First of all, we observe that if a point $P$ in $V_\Pi\cap V_{\Pi'}$ 
does not belong to the scroll $\Sigma$, then it belongs to the fundamental 
locus. Indeed in this case
$$
P\in \ell\cap \ell', \ \text{where}\ \ell\in G_\Pi,\ \ell\in G_{\Pi'},
\ \text{and}\ \ell\neq \ell' 
$$
---as in the proof of Proposition \ref{prop:ai}, $G_\Pi$ and $G_{\Pi'}$ 
denote the subvarieties of the Grassmannian corresponding to the two scrolls
$V_\Pi$ and  $V_{\Pi'}$. Since $\Pi$ and $\Pi'$ are general, 
$P$ belongs to infinitely many lines of $\Lambda$.

Reciprocally, if $P\in \Psi$ is a general focal point, the set of the lines of 
$B$ through $P$, $\chi_P$, is a cone of dimension (at least) two, so its 
intersection with $\Pi$ and $\Pi'$ is not empty and therefore 
$P\in V_\Pi\cap V_{\Pi'}$. 

Finally it remains to prove  that the degree of the scroll is $a_0+2a_1+a_2$ 
if $n>3$; but this follows from the Schubert calculus:
\begin{equation*}
[G_\Pi]\cdot [G_{\Pi'}] = \sum_{i=1}^\nu
((\sum_{j=0}^{i-1}(\nu-i+1)a_j+\sum_{j=i}^{\nu}(\nu-j+1)a_j)
\sigma_{n-i,i+1}).
\end{equation*}
Since the scroll $\Sigma$ has dimension $(n-2)$, to obtain its degree we can 
calculate the intersection of $[G_\Pi]\cdot [G_{\Pi'}]$ with the Schubert cycle
$\sigma_{n-3}$, \ie the lines which meet a general plane, and the result 
follows from the Schubert calculus. 
\end{proof}

Since in the applications of the preceding proposition the cases for $n\le 3$ 
can be treated analogously to the other ones (and the main results are the 
same), we will consider from now on only the cases with $n>3$.


\subsection{The parasitic schemes}
\label{sec:par}
A central definition, introduced first in \cite{DP2}, is that of
\emph{fundamental $d$-locus} for a first order congruence, 
\ie the subscheme of the fundamental locus of 
pure dimension $d$, with $0\le d\le n-2$, which is met by the general line of 
the congruence. Let us see how these schemes are 
constructed: the closed set
$$
S_d:=\{(\Lambda(b),P)\in\Lambda\mid \rk(\ud f_{(\Lambda(b),P)})\le d\}
$$
has a natural subscheme structure, which is defined by a Fitting ideal, \ie
the ideal generated by the $(d+1)$-minors of $\ud f$ (or, by the Fitting 
lemma, see \cite{E}, by the $d$-minors of $\lambda$), see \cite{SK}; in 
particular, $S_{n-1}=V$. Let us define 
$$
D_{d+1}:=\overline{S_{d+1}\setminus S_d}
$$ 
with the scheme structure induced by $S_{d+1}$. Finally, we consider the 
scheme-theoretic image $\Phi_d$ of $D_{d+1}$ in $\p^n$ under $f$. The 
component of $\Phi_d$ of pure dimension $d$ (with the scheme structure induced 
by $\Phi_d$) which is met by the general line of the congruence is the 
fundamental $d$-locus. 
But the fundamental $d$-loci do not always fill up the fundamental locus:
for example, if we consider the congruence in $\p^4$ given by the trisecant 
lines of a surface which contains a plane curve of degree at least three, the 
plane of this curve is contained in the fundamental locus, because all its lines 
are focal, but it is not a component 
of the fundamental $2$-locus. The components we miss are the parasitic schemes:

\begin{definition}\label{df:par}
An irreducible subscheme $\eta$ of the focal locus $\Phi$ of a congruence 
$(\Lambda,B,p)$ of (pure) dimension $d$, with $2\le d\le n-2$, is called 
\emph{$i$-parasitic for the congruence $(\Lambda,B,p)$} (or simply 
\emph{parasitic}) if  
\begin{enumerate}
\item through every point $P\in \eta$, there pass infinitely many (focal) 
lines of the congruence contained in $\eta$,  

\item $\eta$ is a component of $\Phi$ with geometric multiplicity 
$i=\length_{\OO_{(\eta)_{\red},\Phi}}(\OO_{\eta,\Phi})$, 

\item  $\eta$ is not met by the general line of the congruence. 

\end{enumerate}
\end{definition}

\begin{remark}
If $\eta$ is an $i$-parasitic $d$-dimensional scheme, it is straightforward 
that its Fano scheme $\F(\eta)$, \ie the family of lines contained in it, 
is such that $d\le\dim(\F(\eta))\le n-2$. In particular, if $n=4$, the support 
of $\eta$ is a plane, \ie we can have only parasitic planes. 
\end{remark}

Actually, in order to classify the first order congruences, the fundamental $d$-loci are 
the important subschemes of the focal locus, since 
such a congruence can be characterized as the set of lines 
which meet the fundamental $d$-loci a certain number of times, thanks to the 
Classification Theorem~\textup{3.2} of \cite{DP2}. 
For stating it, we need the following notations: we will denote by $C_i$ the 
fundamental $(n-1-i)$-locus, and by $C_{i,j}$ its irreducible components (with 
the scheme structure induced by $C_i$). Moreover, if $\Lambda(b)$ is a  line 
of the congruence not contained in  $C_{i,j}$, we will set  
$C_{i,j}\cap\Lambda(b):= Z_{i,j}=P_1\cup\dotsb\cup P_{k_{i,j}}$, where 
$P_a$, $a\in\{1,\dotsc, k_{i,j}\}$ is a fat point of length $h_a$. Then, the 
Classification Theorem is the following: 

\begin{theorem}\label{thm:relfu} 
Let $(\Lambda,B,p)$ be a first order congruence of lines in $\p^n$. 
Then, using the 
above notations, the congruence $B$ is a subvariety (and so it is 
irreducible and of dimension $n-1$)  
of the set $B_1\subset\Gr(1,n)$ 
of lines 
meeting $C_{i,j}$ in a zero-dimensional scheme $Z_{i,j}$, for every 
$i=1,\dotsc, n-1$, and $j=1,\dotsb,j(i)$ and the following relation holds:

\begin{equation}\label{relfu}
n-1=\sum_{i=1}^{n-1}(\sum_{j=1}^{j(i)} (\sum_{a=1}^{k_{i,j}}h_a)), 
\end{equation}
and $h_a\ge i$. 
\end{theorem}

See \cite{DP2} for a proof. Actually, in \cite{DP2} we made the implicit 
assumption that $B_1$ was irreducible and of dimension $n-1$; 
but the proof given there holds also in the present case.  

Actually, we do not know examples of first 
order congruences $B$ that 
are only a component of $B_1$ (and not the 
whole set), so we conjecture that this is irreducible, if it has order one. 
Instead, if $a_0>1$ (and $\Phi\supsetneq\Psi$) there are examples in which 
$B_1$ is reducible:
\begin{example}
Let us consider a smooth surface $F\subset\p^4$ 
of degree $10$ and sectional 
genus $6$, \ie either an abelian or a bielliptic surface and 
$F$ is fibred in plane cubics and the trisecants of 
the cubics generate a quintic hypersurface (see \cite{rank}) and therefore 
they generate a congruence of order zero. But it is easy to see that there 
is (at least) another component of $B_1$: in fact, if we apply the 
triple point formula to $F$ (see for example \cite{Au}, 
\ie Formula~\eqref{eq:au} which can be found below), we obtain that 
$B_1$ has to have order $40$.  

We have $\Phi\supsetneq\Psi$ since if we apply the quadruple point formula for 
a smooth surface in $\p^4$ (see \cite{Bar}) to $F$, 
we obtain that the degree of the hypersurface $Q$ of the quadrisecant lines 
of $F$ is $130$ (and clearly $Q\subset\Phi$, since a quadrisecant is a focal
line). 

We could conclude that $\Phi\supsetneq\Psi$ also applying 
Theorem~\ref{thm:dgb} below. A similar discussion about this example can be 
found in \cite{ADR}. 
\end{example}

In any case, the results of \cite{DP2} are correct, since there the 
congruences are given by the joining lines a curve and a surface (and then 
$B$ is birationally equivalent to the product of the curve and the surface).

\begin{proposition}\label{cor:fu}
The components of the focal locus of a first order congruence 
of lines in $\p^n$ are the fundamental $d$-loci, where $d$ varies between $0$ and 
$n-2$, and the parasitic schemes. 
\end{proposition}

\begin{proof}
If $\eta$ is a (maximal) subscheme of the focal locus of dimension $d$ which 
is  not a fundamental $d$-locus and $P$ a general point in it, then through 
$P$ there pass infinitely many lines of the congruence. Besides, since $P$ is 
not a point in the fundamental $d$-locus, on each of the infinitely many lines 
of the congruence through $P$ there are at least $n-1$ other foci distinct 
from $P$; therefore, these are fundamental lines, which are to be contained in 
$\eta$. 
\end{proof}

\begin{remark}
The definition of parasitic plane for a congruence of lines in $\p^4$ was 
introduced by Ascione in \cite{A}, and used in \cite{M1}. Ascione was 
interested in classifying the surfaces $S$ with one apparent triple point, \ie 
the surfaces such that for a general point in $\p^4$ there passes only one 
trisecant line of $S$. Then the parasitic planes come out naturally as the 
components of the focal locus different from $S$.

Marletta generalized this concept to general congruences in higher dimensional 
spaces, in \cite{M2}: he defined an $i$-parasitic $d$-plane, with $d<(n+1)/2$ 
as a linear space such that a general line $\ell$ in it is of multiplicity $i$ 
for the congruence (\ie  $\deg(f^{-1}(\ell))= i$). The limitation on the 
dimension is due to the fact that higher dimensional linear spaces contain 
``too many'' lines, \ie they contain families of lines of dimension greater 
than $n-1$. 
\end{remark}

\begin{remark} 
An $i$-parasitic $(n-2)$-dimensional scheme $\eta$, by definition, is a 
component of the focal locus $\Phi$ whose geometric multiplicity in $\Phi$ is 
$i$. Then it is a component of the scheme $\Gamma$ of Proposition~\ref{thm:pp} 
with geometric multiplicity $i^2$. In fact, the intersection multiplicity 
$i((\eta)_{\red},V_\Pi\cdot V_{\Pi'},\p^n)$ of $(\eta)_{\red}$ in 
$V_\Pi\cdot V_{\Pi'}$ 
is equal to the geometric multiplicity of  $(\eta)_{\red}$ in $\Gamma$, 
but $i((\eta)_{\red},V_\Pi\cdot V_{\Pi'},\p^n)=i^2$. 
\end{remark}

From now on, we will denote by $i$ the multiplicity of a general parasitic 
$(n-2)$-dimensional scheme. We will also set $x:=\sum i^2b_i$, where $i$ 
varies among all the $i$-parasitic $(n-2)$-dimensional schemes $\eta_i$ 
and $b_i:=\deg(\eta_i)_{\red}$. 

\begin{definition} The union  of the fundamental $d$-loci of  $\Psi$ is called 
\emph{pure fundamental locus}, or, in what follows, simply 
\emph{fundamental locus} and it is denoted by $F$. 
\end{definition}

\begin{proposition}\label{prop:pas}
If $F$ is the pure fundamental locus and $\eta$ is an $i$-parasitic 
$d$-dimensional scheme, then $F\cap \eta$ has codimension one in $\eta$. 
Besides, 
\begin{equation}\label{eq:binomg}
i=\sum_{k=1}^{n-2}\sum_{j=1}^{j(k)}\binom{\mu_{k,j}}{n-1}
\end{equation}
where $\mu_{k,j}:=\deg(C_{k,j}\cap\eta)$ and $C_{k,j}$, $j=1,\dotsc, j(k)$, 
varies among the irreducible components of the fundamental $k$-locus.
\end{proposition}

\begin{proof}
If, by contradiction, $\dim(F\cap \eta)\le d-2$ and  $\ell$ is a general line 
contained in $\eta$, we have $\ell\cap F=\emptyset$. But every focal line 
contained in the scheme $\eta$ is a line of the congruence, so it contains at 
least $(n-1)$ focal points. 

A general focal line $\ell$ contained in the parasitic scheme $\eta$ 
intersects $C_{k,j}$ in $\mu_{k,j}$ points, so taking the $\mu_{k,j}$ points 
$(n-1)$ by $(n-1)$ we obtain the formula.
\end{proof}


\subsection{General formulae and degree bounds}
\label{sec:db}
We will suppose, throughout this subsection, that the fundamental 
$(n-2)$-locus $C\subset F$ is not empty. 

In what follows, let $C_j$, $j=1,\dotsc,h$ be the $(n-2)$-dimensional 
irreducible components of $C$; we will denote then by $m_j$ the degree of 
$(C_j)_{\red}$ and by $k_j$ the algebraic multiplicity of $(C_j)_{\red}$ on 
$V_\Pi$. 
Finally, we put $\ell_j:=\length((C_j)_{\red}\cap\Lambda(b))$, where 
$\Lambda(b)$ is a general line of the congruence. 

\begin{proposition}\label{prop:agen}
The following formulae hold:
\begin{align}
\sum_{j=1}^h\ell_jk_j& \le  a_0+a_1,\label{eq:agen}\\
(a_0+a_1)^2&= x+\sum_{j=1}^hk^2_jm_j+a_0+2a_1+a_2; 
\label{eq:bgen}
\end{align}
moreover, if the pure fundamental locus has pure dimension $n-2$, the equality holds 
in formula~\textup{\eqref{eq:agen}}.
\end{proposition}

\begin{proof}
If we take a line $\Lambda(b)$ of the 
congruence not contained 
in $F\cap V_\Pi$, then, intersecting $\Lambda(b)$ with $V_\Pi$, we obtain a 
zero dimensional scheme of length $a_0+a_1$, since, as 
we have seen in Proposition~\ref{prop:ai}, this is the degree of $V_\Pi$. 
$V_\Pi\cap \Lambda(b)$ contains the schemes $(C_j)_{\red} \cap\Lambda(b)$, 
which have length $\ell_j$, and the intersection multiplicity of 
each of the irreducible components of 
$(C_j)_{\red} \cap\Lambda(b)$ in it is $k_j$, 
so the relation~\eqref{eq:agen} is proved. If $F$ has pure dimension $n-2$, 
then $\Lambda(b)$ does not intersect 
$V_\Pi$ in points different from $F$, and the equality in 
formula~\textup{\eqref{eq:agen}} holds. 

We recall that, by Proposition~\ref{thm:pp}, 
the degree of $\Gamma:=V_\Pi\cap V_{\Pi'}$ is $(a_0+a_1)^2$, 
and its components of maximal dimension 
are---by Proposition~\ref{cor:fu}---the 
fundamental $(n-2)$-locus, the parasitic schemes and the scroll $\Sigma$. 
A component $(C_j)_{\red}$ of the fundamental $(n-2)$-locus has geometric 
multiplicity in $\Gamma$ equal to $k_j^2$. The proof of this fact is 
as the one of the preceding remark for the multiplicity of the 
$i$-parasitic schemes: 
the intersection multiplicity 
$i((C_j)_{\red},V_\Pi\cdot V_{\Pi'},\p^n)$ of $(C_j)_{\red}$ in 
$V_\Pi\cdot V_{\Pi'}$ is equal to the geometric 
multiplicity of  $(C_j)_{\red}$ in $\Gamma$, 
but $i((C_j)_{\red},V_\Pi\cdot V_{\Pi'},\p^n)=k_j^2$. 

Finally, as we have seen, the $i$-parasitic schemes of dimension  $n-2$ 
give $x=\sum i^2b_i$, and the scroll $\Sigma$ has degree 
$a_0+2a_1+a_2$, so we get formula~\eqref{eq:bgen}.
\end{proof}

\begin{theorem}\label{thm:dgb}
If $(\Lambda,B,p)$ is a first order congruence of lines in $\p^n$ 
such that its pure 
fundamental locus $F$ is irreducible and 
coincides with the fundamental $(n-2)$-locus; then 
\begin{equation*}
\frac{n-1}{k'}<m<(n-1)^2, 
\end{equation*}
where $m:=\deg (F)_{\red}$ and $k'$ is the geometric multiplicity of
$(F)_{\red}$ in $F$. 
\end{theorem}

\begin{proof}
If we substitute formula~\eqref{eq:agen}, which in this particular case becomes 
an equality, 
in formula \eqref{eq:bgen}, we obtain
\begin{equation}
(n-1)^2k^2 -k^2m-1-2a_1-a_2= x\ge 0,\label{eq:a2gen}
\end{equation}
---where $k$ is the algebraic multiplicity
of $(F)_{\red}$ on $V_\Pi$---and since  $-1-2a_1-a_2< 0$ by formula 
\eqref{eq:bgen}, we deduce $m<(n-1)^2$. 

We have that $n-1<k'm$ by degree reasons, since the congruence is given by the 
lines which intersect $F$ in a zero dimensional scheme of length  
$n-1$. 
\end{proof}

An interesting corollary of the preceding theorem is the following (see Theorem~1.3 of 
\cite{DP3}): 

\begin{corollary}
Given a first order congruence in $\p^n$, which is moreover generated by the 
$(n-1)$-secant lines of an irreducible $(n-2)$-dimensional variety $F$ of degree$d$, 
then $F$ cannot be a complete intersection. 
In particular, if $F$ is smooth, $n\le 5$. 
\end{corollary}

See Theorem 1.3 of \cite{DP3} for a proof. 




\section{Congruences in $\p^4$}\label{sec:1t}
\subsection{Multiple point formulae}
We start with some multiple point formulae for a smooth 
surface $S\subset\p^n$. In this subsection, let us denote by 
$m$ the degree of $S$, by $\pi$ its sectional genus and by $\chi$ 
its Euler-Poincar\'e characteristic. Besides, $K$ and $H$ are the canonical 
and hyperplane divisors, respectively. A \emph{trisecant} $(n-3)$-plane  
of $S$ is a linear space that intersect the surface in a zero dimensional 
scheme of length (at least) three. 
A standard application of the triple point formula of \cite{K1} is the 
following: 

\begin{proposition}\label{prop:3ple}
Let $S$ be a smooth surface in $\p^n$ and $\ell$ be an 
$(n-4)$-plane with $\ell\cap S=\emptyset$ and   
through which there is a finite number of trisecant 
$(n-3)$-planes of $S$. If we set  
$t:=\length(\xi\cap W_P)$ where 
$\xi, W_P\subset\Gr(n-3,n) \subset \p^{\binom{n+1}{3}-1}$ 
are, respectively, the family of trisecant $(n-3)$-planes of $X$ and the 
Schubert variety of the $(n-3)$-planes through $\ell$ 
(embedded via the Pl\"ucker 
embedding), then the following formula holds:
\begin{equation}
t=\frac{1}{6}m^3-\frac{3}{2}m^2+\frac{13}{3}m-m\pi
+K^2+8\pi-4\chi-8.
\label{eq:3ple}
\end{equation}
In particular, for the general $\ell$, there will pass  
$t$ trisecant $(n-3)$-planes.
\end{proposition}

The proof is analogous to the one given for Proposition 3.1 of \cite{DP3}. 
Clearly, $t$ is the number of triple points of the general 
projection of $S$ from 
$\ell$ to a $\p^3$, \ie $t$ is the order of the congruence of the trisecant 
lines of a general projection of $S$ in $\p^4$. 
We note that in the case $n=4$ we obtain 
\begin{equation}\label{eq:au}
t=\binom{m-1}{3}-\pi(m-3)+2\chi -2,
\end{equation}
which can be found in \cite{Au}. 



\subsection{Congruences}
In this section we study the congruences of lines of order one with a 
reduced irreducible surface $F$, 
proving Theorems~\ref{thm:primo} and \ref{thm:smo}. 


In this case the sequence of degrees is formed only by the order, which is
one, and the class. In what follows, we will use the following notations: 
\begin{itemize}
\item $F\subset\p^4$ is an irreducible surface of degree $m$ and 
sectional geometric genus $\pi$; so $3<m<9$ by Theorem \ref{thm:dgb}; 
\item $\alpha:Y\rightarrow F$ the normalization of $F$ and 
$\beta:S\rightarrow Y$ the minimal desingularization of $Y$. We also set 
$\varphi:=\beta\circ\alpha$; 
\item $K$ denotes the canonical divisor of $S$ and $\chi$ 
its Euler-Poincar\'e characteristic; besides, we set $H:=\varphi^*(\OO_F(1))$, 
and therefore $2\pi-2=m+H\cdot K$. 
\item $H'$ is a general hyperplane in $\p^4$ and  
$C_{H'}:=F\cap H'$ is the general hyperplane section of $F$; 
\item $h$ is the number of apparent double points of the curve $C_{H'}$.  
\item $\eta$ is a general $i$-parasitic plane, 
$C=\eta\cap F$ is the curve in $\eta$ of Proposition~\ref{prop:pas}; 
$\mu\ge 3$ is the degree of $C$, so $i=\binom{\mu}{3}$, 
by Proposition~\ref{prop:pas}; 
\item $x:=\sum i^2$ where the sum varies among all the 
parasitic planes of the congruence, $a$ is the class of the congruence 
and $k$ is the algebraic multiplicity
of $F$ on $V_\Pi$ (see Proposition~\ref{prop:ai}); so, 
formula~\eqref{eq:a2gen} becomes 
\begin{equation}\label{eq:a2}
(n-1)^2k^2 -k^2m-1-2a= x\ge 0. 
\end{equation}
\end{itemize}

We begin with some general results about congruences and trisecant lines of 
surfaces in $\p^4$. 

\begin{proposition}\label{ex:trisec} 
The family $B$ of the trisecant lines of a nondegenerate surface 
$F\subset\p^4$ is either empty (with the convention that a line contained 
in $F$ is not a trisecant) or has dimension three. 
\end{proposition}

\begin{proof}
There is a natural surjective map $\phi:\Al^3(F)\rightarrow B$, where 
$\Al^3(F)\subset \Hilb^3(F)$ is the subscheme of collinear 
(\ie which are contained in a line) triple points (see \cite{LEB}).

Indeed, $\Al^3(F)\cong \Hilb^3(F)\times_{\Hilb^3(\p^4)}\Al^3(\p^4)$, so 
either $\Al^3(F)$ is empty or has dimension at least three, because 
$\dim(\Al^3\p^4)=9$. Hence the same conclusion holds for $B$, because the 
general fibre of $\phi$ is finite (if not, all the trisecants lines of $F$ 
are contained in $F$, and $F$ is a plane). 
 
Moreover, we can exclude that $\dim(B)>3$. 
In fact, if it had dimension four, then through a general point $P$ of
$F$ would pass $\infty^2$ lines of the congruence. Therefore, 
this would be the join of $S$ and $P$ and so the general secant 
line would be a trisecant: 
this would contradict the trisecant lemma. 
\end{proof}

We need also the following 
lemma, which was suggested by K. Ranestad. 

\begin{lemma}\label{lem:cone}
If $F\subset\p^n$ is an $(n-2)$-dimensional irreducible cone with vertex 
$V\cong\p^k$, then the family of $(n-1)$-secant lines of 
$F$ either is a congruence of order zero or it has dimension greater than 
$n-1$. 
\end{lemma}

\begin{proof}
If $P\in\p^n$ is a point through which there is an $(n-1)$-secant line 
$\ell_P$ of $F$, then every line of $\overline{\ell_P V}$ passing through $P$ 
is in fact an $(n-1)$-secant line of $F$. 
\end{proof}

We give now some remarks about the congruences generated by the trisecants 
of a singular surface in $\p^4$. 

\begin{definition}
Let $F$ be an irreducible surface in $\p^n$; the number of its \emph{apparent 
triple points}---which will be denoted in what follows by $\tau$---is the number 
of triple points of a general projection of $F$ to a $\p^3$ which are not 
projections of triple points of $F$.  
\end{definition}

Therefore, the number of apparent triple points $\tau$ is equal to the $t$ of 
Proposition~\ref{prop:3ple} for a smooth surface 
$S\subset\p^n$.

\begin{proposition}
If $F$ is an irreducible surface in $\p^4$, then $\tau\ge a_0$---where $a_0$ is 
the order of the congruence of the trisecants of $F$; moreover, 
if $F$ has isolated singularities, the equality holds.   
\end{proposition}

\begin{proof}
It follows from the fact that a projection of a 
trisecant line gives a triple point, and 
if a projeciton of a secant line---which is not a trisecant---gives a 
triple point, then one of the two secancy point is a multiple point for $F$ 
and so lies on a singular locus of dimension one.
\end{proof}

The two numbers $a_0$ and $\tau$ are not always the same, as the following 
example shows: 

\begin{example}\label{ex:1}
It is a well known fact, which also follows easily from 
formula~\eqref{eq:3ple} that the Veronese surface $V\subset\p^5$ has one 
apparent triple point. In fact, the projection of $V$ from a line that 
does not intersect it to a $\p^3$ is the Steiner Roman surface. 
Therefore, if we project it from a point not in the 
secant variety of $V$ to a $\p^4$, we obtain a smooth surface $F$ whose 
trisecants generate a first order congruence. 

But if the projection is not general, \ie if we project from a point 
(not on the Veronese surface) of the secant variety of $V$, 
the projection $F\subset \p^4$ has not 
trisecants, 
since it is easy to see that $F$ is the complete intersection of two 
quadrics. Moreover, the singular locus of $F$ is a (double) line. 
\end{example}

\begin{corollary}\label{cor:a0t}
Let $S$ be a smooth surface in $\p^n$, with $n\ge 5$; for a particular proiection $F$ of $S$ 
from an $(n-5)$-plane $\ell$ with $\ell\cap S=\emptyset$ to a $\p^4$, 
the order $a_0$ of the 
family of the trisecant lines of $F$ is less than or equal to the number of 
apparent triple points $\tau$ of $S$---with $\tau=a_0$ if $F$ has only isolated 
singularities.  
\end{corollary}

\begin{proposition}
Let $F$ be an irreducible surface in $\p^4$ with $\sing(F)=:\ell$ of dimension 
one and whose trisecant lines generate a first order congruence 
$B$. Then we cannot have 
$\length(\lambda(b)\cap \ell)=2$, with $\lambda(b)\in B$   (\ie two of the 
foci contained in $\ell$ as a subscheme in $F$).  
\end{proposition}

\begin{proof}
In fact, if it were so,(an irreducible component of) $B$ would be the join of 
(a component of) $\ell$ and $F$; 
$(\ell)_{\red}$ 
must be a plane curve, since its secant variety is contained in the focal 
locus. Let $\pi$ be a plane such that $(\ell)_{\red}\subset\pi$; 
in each element $\Pi$ of the pencil of hyperplanes $\p^1_\pi$ 
containing it, the congruence induces a first order congruence of 
$\p^3(\cong\Pi)$. But this cannot happen, since $\Pi\cap F$ contains 
$\ell$ at least as a double curve and therefore by the classification of 
first order congruences in $\p^3$ of 
\cite{DP1},  $(\Pi\cap F)_{\red}=(\ell)_{\red}$, \ie $(F)_{\red}$ is a plane. 
\end{proof}

\begin{proposition}
The family of the tangent lines to a nondegenerate surface $F\subset\p^4$ 
cannot have $F$ as its sole component of the (pure) fundamental locus. 
\end{proposition}

\begin{proof}
If it were so, since on every line $\ell_P$ which is 
tangent to a $P\in F$, we would 
have three foci, there would be another point $Q\in F\cap \ell_P$, $Q\neq P$, 

and therefore every tangent plane $\T_P F$ 
would intersect $F$ in a (plane) curve $C_P$---which is necessarily singular 
in $P$. 
Therefore, we would have a family of dimension two of plane curves, and 
by C. Segre's Theorem (Theorem 4 of \cite{MP}), this imply that $F$ is the 
Veronese surface, the cubic scroll or a cone. But none of these cases can 
occur since---for example---$C_P$ should always singular in $P$.   
\end{proof}



\begin{theorem}\label{thm:6}
If the trisecant lines of a surface $F$ in $\p^4$ generate a (possibly 
reducible) congruence of lines of order $a_0>1$, then $\deg(F)=m\ge 6$. 
\end{theorem}

\begin{proof}
If $P\in\p^4$ is a general point, there will pass at least 
two trisecant lines through it, $\ell_1,\ell_2$, and the plane $\Pi$ generated 
by them intersects $F$ in at least six points. Therefore, if $m\le 5$, 
$C_P:=\Pi\cap F$ is a curve. In such a way, we obtain a family of 
dimension (at least) two of planes, and each of them containing a plane curve 
of $F$, \ie a congruence of planes in $\p^4$ (see \cite{SC}). 
Clearly, $\deg(C_P)<3$, since otherwise every trisecant is a secant, 
which contradicts the trisecant lemma (or Proposition~\ref{ex:trisec}). If 
$\deg(C_P)=2$, then $F$ is a rational variety of degree at least $5$, since 
$\Pi$ contains at least two points in $F$ out of the conic. 
But $F$ contains a family of dimension two of conics, so by C. Segre's 
Theorem, 
$F$ can be the Veronese surface or a cubic scroll 
or a cone and by degree reasons and by Lemma \ref{lem:cone} it is easy to show 
that none of these cases can occur. 

If $\deg(C_P)=1$ again by C. Segre's Theorem we have that, given one of 
these lines $\ell\subset C_P$, there will pass infinitely many planes 
$\Pi$'s through it, \ie $\ell\subset C_P$ for infinitely many $P$'s. 
If there is only finitely  many of these lines, then there exists only one: 
the family of planes containing a line has in fact dimension two. 
But then $\ell$ is a component of the fundamental $1$-locus, 
which cannot occur, since otherwise we would have more than three foci on 
each line of the congruence. 
Therefore $F$ is a ruled surface. It is rational also: fix a general line $r$ 
of $\p^4$ and consider the map which associates to $P\in r$ the corresponding 
line $\ell\subset C_P$. So $F$ is (the projection of) a rational normal 
scroll. Then, applying formula \eqref{eq:3ple} (and Corollary~\ref{cor:a0t}) 
to the rational normal scrolls we get the result. 
\end{proof}


We need the following lemma, and for proving it we need some results on joins 
of projective varieties, which can be found in \cite{AA}. We will follow the 
multiplicative notation introduced there: if $X,Y\subset\p^n$ the join of 
$X$ and $Y$ will be denoted by $XY$. 

\begin{lemma}\label{lem:ad}
Let $C_{H'}$ be a general hyperplane section of the fundamental surface $F$ and 
$P$ a point in it. Then the tangent space $\T_P C_{H'}$ has dimension at most two.
\end{lemma}

\begin{proof}
If there were a point $P$ for which $\dim(\T_P C_{H'})=3$, $F$ would contain a 
singular curve $D$ (with $P\in D$) such that if $Q\in D$, then 
$\dim(\T_P F)=4$. The join $DF$ of $D$ and $F$ cannot be the whole $\p^4$, 
otherwise $D$ would be a fundamental $1$-locus. So $\dim DF=3$, therefore 
$DF$ is a cone with vertex $D$. In particular, $D$ is a line and $DF$ has, 
as basis, a curve $E$. If $Q$ is a point in $D$, 
we must have  $\dim(QE^2)=2$ and so $E$ must be a plane curve 
and $Q$ has to be contained in the plane $E^2$ ($E$ cannot be a line), 
which is absurd, since we would have $D\subset E^2$. 
\end{proof}

From now on, we will assume that the family of the trisecant lines of 
$F$ is irreducible.

\begin{proposition}
Let $C_{H'}$ be a general hyperplane section of the fundamental surface $F$. 
Then, if $h$ is the number of the apparent double
points of $C_{H'}$, we have
\begin{align}
a&=h(m-2)-\binom{m}{3},\label{eq:zeuthen}\\ 
k&=h-m+2,\label{eq:c}
\end{align}
and if $m\neq 5$, we obtain:
\begin{equation}\label{eq:e}
h=\frac{m(m+2)}{6}-1.
\end{equation}
\end{proposition}

\begin{proof}
A formula due to Cayley---and to Le Barz, in the singular case---(for a proof 
see \cite{LB}) which expresses the number of trisecant lines 
(\ie lines containing 
zero dimensional subschemes of length three contained in $C_{H'}$) 
of a curve in $\p^3$ meeting a line, in our situation gives---thank to 
Lemma \ref{lem:ad}---formula~\eqref{eq:zeuthen}. 

Since $k$ is the algebraic multiplicity of $F$ on $V_\Pi$, 
this means that through a 
general point $P$ in $F$ there are $k$ lines of the congruence (\ie trisecant 
lines of $F$) 
that meet $\Pi$, and these lines belong to the $\p^3$ spanned by 
the point $P$ and the plane $\Pi$. 

Formula~\eqref{eq:c} can be obtained by computing the (geometric) 
genus of the generic hyperplane section $C_{H'}$ in two ways
with the Clebsch formula: 
the first by projecting $C_{H'}$ to a plane from a point
not belonging to it, and the second by projecting from $P\in F$ 
to a general plane. 

Formula~\eqref{eq:e} is obtained from 
equalities~\eqref{eq:agen},~\eqref{eq:zeuthen} and 
\eqref{eq:c} eliminating $a$ and $k$, and simplifying by $m-5$. 
\end{proof}





\begin{corollary}\label{cor:456}
The possible values for the degree of the fundamental surface $F$ are  
$m=4,5,6$.
\end{corollary}

\begin{proof}
This is due to the fact that, by Theorem~\ref{thm:dgb}, $3<m<9$ and 
for the values $m=7,8$ we have non-integer 
values for $h$ in formula~\eqref{eq:e}.
\end{proof}

So, the only possible cases to analyse are $m=4,5,6$.


\subsection{Case $m=4$} 
We will prove (in Theorem~\ref{thm:ver}, which in fact precises 
point~\eqref{ve} of Theorem~\ref{thm:primo}) that the only 
irreducible surface in $\p^4$ whose trisecant lines
generate a first order congruence is
a smooth projected Veronese surface.




A standard fact---whose proof can be found in \cite{Fuj}---is the following

\begin{lemma}\label{cor:kc}
Let $S\subset\p^N$ be an irreducible surface 
with rational sections; then $S$ is a projection 
of either a rational normal scroll or a Veronese surface.
\end{lemma}

From this we can classify the first order congruences whose fundamental
surface has degree four:

\begin{theorem}\label{thm:ver}
The only surface $F$ of degree four 
which generates a first order congruence is the projection of the 
Veronese surface $S$ from a point not lying on the secant variety of $S$. 
In particular, $F$ is smooth. 

\textit{Vice versa}, the trisecant lines of the smooth, projected Veronese 
surface in $\p^4$ generate a first order congruence.
\end{theorem}

\begin{proof}
By formula~\eqref{eq:e} and by the Clebsch formula, $F$ has rational sections; 
therefore, by Lemma~\ref{cor:kc}, the possibilities are the projection of the 
Veronese surface, which has in fact degree four, and three quartic scrolls,
projection of the rational normal scrolls of degree four in $\p^5$. By 
Lemma \ref{lem:cone} the conic scroll cannot occur. 

Then, we apply Example \ref{ex:1} to get that the trisecants of the 
smooth Veronese surface generate a first order congruence (and the 
trisecants of the singular one do not), while the trisecants of 
the scrolls generate a congruence of order zero: in fact, we have 
$\chi=1$ for all these surface, while for the Veronese surface we have $K^2=9$ 
and $K^2=8$ for the scrolls.  

There are not (non-secant) lines in $\p^5$ of the Veronese surface $S$ 
for which there are infinitely many trisecant planes by degree reasons. 
%
%
\end{proof}

\begin{remark}
Besides, we have also that, for the Veronese surface,
$a=2$, $k=1$ and $x=0$, 
and so we do not have parasitic planes.

The fact that the scrolls do not generate a first order congruence 
can be proven with the following nice geometric argument, suggested to us 
by Andrea Bruno: if we project a scroll $S\subset\p^5$ of degree four (\ie all 
the surfaces of minimal degree but the Veronese surface)
from two general points $P_1$ and $P_2$ lying on a secant line $\ell$ in $S$
(or, which is the same, from the secant line $\ell$ given by the
span of the two points), we obtain a quadric $Q$ in
$\p^3$. So, the projection of $S$ from one point $P_1$ is
contained in the quadric cone given by the quadric $Q$ and the other
point $P_2$. Since any point in $\p^5$ lies on a secant line of the 
scroll, the previous projection is general. Then the scrolls
cannot generate a first order 
congruence since a trisecant line of $F$ is contained in the quadric cone.

\end{remark}

\subsection{Case $m=5$}
In this subsection we will prove Theorem~\ref{thm:m5}, which is a more precise 
version of points~\eqref{dp} and \eqref{sc} of Theorem~\ref{thm:primo}. 


\begin{theorem}\label{thm:m5}
If $F$ is an irreducible surface in $\p^4$ of degree five whose 
trisecant lines generate a first order congruence, then 
\begin{enumerate}
\item\label{m5:1} 
either $F$ has sectional genus one, in which case it 
is the projection of a---possibly singular---Del Pezzo surface $D$ 
of $\p^5$ (see 
for example \textup{\cite{dem}}) from a point neither contained in $D$ nor in a 
plane containing an irreducible conic of $D$, 
and we have $5$ $1$-parasitic planes, 
\item\label{m5:2} or $F$ has sectional genus zero and it is a projection of 
a smooth rational normal scroll $S$ in $\p^6$  
from a non-secant line 
of 
$S$. 
Moreover we have one $4$-parasitic and three $1$-parasitic planes. 
\end{enumerate}
\textit{Vice versa}, the trisecants of such surfaces generate a first order 
congruence.
\end{theorem}

\begin{proof}
First of all we observe that $F$ cannot be smooth.
In fact, if $F$ were smooth, the number of trisecant lines $t$ 
of $F$ through a general point in $\p^4$ would be
$$
t=2(\chi-\pi+1),
$$ 
by formula~\eqref{eq:au}, so $t$ is an even number and then it cannot be one.

Then, we deduce that we have only two possible lists of invariants, 
which are the following 
\begin{align} 
h&=5, & k&=2, & a&=5, & x&=5,  &\pi&=1, \label{eq:h5}\\ 
h&=6, & k&=3, & a&=8, & x&=19, &\pi&=0.   \label{eq:h6}
\end{align}
In fact, from \eqref{eq:a2}, solving the inequality, we have that $k>1$, 
and from \eqref{eq:c} we deduce that $k=h-3$, and finally that $h>4$.
On the other hand, by the Clebsch formula, $h\le 6$.
From this we deduce \eqref{eq:h5} and \eqref{eq:h6}.

Let us start with case~\eqref{eq:h5}:  
first of all, from \eqref{eq:binomg} and from $x=5$ we deduce that the only 
possibility is to have five $1$-parasitic planes
such that each of them contains 
a cubic plane curve.

From  Theorem 4 of \cite{Hala1} and Lemma \ref{lem:cone}, we get that the 
only possibility is that $F$ is the projection of a (possibly singular) 
Del Pezzo surface in $\p^5$. 



It remains to prove that this surface generates a first order congruence. 
Since $F$ is a Del Pezzo surface,
which is rational, then
$\chi=1$, $K=-H$,  $K^2=5$ and $\pi=1$, and so, by formula \eqref{eq:3ple} 
$t=1$ (If the Del Pezzo is singular, we can take its desingularization and 
apply Kleiman's formula to it: the invariants are the same since it is always 
the plane blown-up in four points). 

We have that the point of projection cannot lie on a plane $\pi_C$ 
spanned by an irreducible conic $C$ contained in the Del Pezzo $D\subset\p^5$ 
since---for example---the 
(only, by Corollary~\ref{cor:a0t} and Theorem~\ref{thm:6}) triple point of 
a projection of $D$ to a $\p^3$ must lie on the projection of $C$ (in fact the 
double locus of the projection has degree $5$ and is connected). 


Let us consider now case~\eqref{eq:h6}. 
By Lemma~\ref{cor:kc}, 
$F$ is a projection $p_\ell$ from a line $\ell$ of a rational normal 
quintic scroll $S$ in $\p^6$. $S$ can be the 
cone $S_{0,5}$ or one of the two smooth scrolls $S_{1,4}$ and $S_{2,3}$;  
but the case of the cone cannot occur by Lemma~\ref{lem:cone}.
We see now that the other two cases occur. 

First of all, we claim that $F$ has one $4$-parasitic and three 
$1$-parasitic planes. Indeed, the general hyperplane section 
$C_{H'}$ of $F$ is a smooth 
rational curve of degree five in $H\cong\p^3$; so, by a 
formula of \cite{LB}, the number of quadrisecant lines
of $C_{H'}$ is one.
Then, from \eqref{eq:binomg} and $x=19$ we obtain that either $F$ 
contains a quartic plane curve and three plane cubics, or $19$ plane cubics.
But every hyperplane section contains a $4$-secant, 
so we are in the first case: indeed, if the quadrisecant lines
were not contained in a plane,  
they would generate a family of dimension two and all these lines would
be contained in the fundamental surface; but the only surface with this 
characteristic is the plane. Therefore the claim is proved.

To conclude it remains to prove that these scrolls generate a 
first order congruence, and by formula \eqref{eq:3ple} we have $t=1$. 

We note that, as in the case of the Veronese, by degree reasons we can 
apply formula \eqref{eq:3ple} for every non-secant $\ell$ in any of the 
cases of the theorem. 
\end{proof}

\begin{remark}
In \cite{PDP} we classified case \eqref{m5:2} of Theorem \ref{thm:m5} 
with the hypothesis that $F$ has isolated singularities and we 
obtained---with a very long geometric argument---that in this case the 
projection must be general; besides, we also obtained a precise description 
of (the mutual position of) the parasitic planes. 
\end{remark}


\subsection{Case $m=6$}\label{eq:6}
We will prove Theorem~\ref{thm:bs}, which precises point~\eqref{bo} of 
Theorem~\ref{thm:primo}. We start with the following

\begin{lemma}\label{lem:bor}
Let $F\subset\p^4$ be an irreducible surface given by the degeneracy locus of 
a $3\times 4$ matrix $M$ of linear forms \ie a \emph{Bordiga surface}.  
Then with some invertible row and column operations, $M$ becomes a matrix 
of the form 
$$ 
M':=\begin{pmatrix}
\ell_{11}&\ell_{12}&\ell_{13}&\ell_{14}\\
\ell_{21}&\ell_{22}&\ell_{23}&\ell_{24}\\
\ell_{31}&\ell_{32}&\ell_{33}&0
\end{pmatrix}, 
$$ 
where $\ell_{i,j}\in \mathbb{C}[x_0,\dotsc,x_4]_1$, 
\ie $M$ is not \emph{$1$-generic} (see \textup{\cite{Ei}}). Moreover, $F$ has at 
most $5$ singular points. 
\end{lemma}

\begin{proof}
In fact, the $2\times 2$ 
minors of the generic $3\times 4$ matrix give the Segre variety 
$\Sigma_{2,3}=\p^2\times \p^3$. If $\ell\cong \p^4$ is such that  
$\dim(\Sigma_{2,3}\cap \ell)=0$, then $\length(\Sigma_{2,3}\cap \ell)\le 5$ 
since these points are on $\Sigma_{2,3}$ and therefore are in 
general position.   
So, applying Proposition~4.1 of \cite{Ei} we get that $M$ cannot be 
$1$-generic.

The $2\times 2$ minors of $M$ cannot determine a curve $r$, since  
the $3\times 3$ minors would contain, 
by---for example---Lemma~2.5 of \cite{hsv} the secant variety of $r$.  
Therefore $r$ could be only a line. But if $r$ were a line, then the 
entries of the matrix $M$ would not generate the vector space 
$\mathbb{C}[x_0,\dotsc,x_4]_1$ (simply evaluate $M$ on $r$) and therefore 
$F$ would be a cone. 
\end{proof}

\begin{theorem}\label{thm:bs}
The only surface $F\subset\p^4$ of degree six whose trisecants generate a 
first order congruence is the Bordiga 
surface in $\p^4$, \ie the degeneracy locus of a map
$\phi\in\Hom(\OO_{\p^4}^{\oplus 3},\OO_{\p^4}^{\oplus 4}(1))$---if 
the minors vanish in the expected (irreducible) codimension two. 
This surface has at most $5$ singular points.

\textit{Vice versa}, the trisecant lines of a Bordiga surface generate a first 
order congruence.
\end{theorem}

\begin{proof}
First of all we can deduce from the results of the beginning of the 
section that if $F$ has degree six,  we have 
\begin{align*}
h&=7, & k&=3, & a&=8, & x&=10, &\pi&=3;
\end{align*}
so $F$ has $10$ $1$-parasitic planes. 
If $\eta$ is one of these planes (and $C=\eta\cap F$ the corresponding cubic 
curve), the general element $\Pi$ 
of the  pencil $\p^1_\eta$ of hyperplanes containing $\eta$ will intersect 
$F$ out of $C$ in a cubic $C_\Pi$ (\ie $\Pi\cap F = C\cup C_\Pi$). 
If $C_\Pi$ were reducible, $F$ would be a rational scroll, but this cannot 
happen, since---for example---$\pi=3$. 

Therefore $C_\Pi$ is a twisted cubic. We can choose coordinates 
$x_0,\dotsc,x_4$ on $\p^4$ such that $\eta=V(x_3,x_4)$ and so  
$\Pi=V(\lambda x_3+\mu x_4)$ where $(\lambda:\mu)\in \p^1$. 
From this we obtain that $C_\Pi$ is the determinantal variety given 
by the $2\times 2$  minors of the following matrix
$$
A_\Pi:=\begin{pmatrix}
\ell_{11}&\ell_{12}&\ell_{13}\\
\ell_{21}&\ell_{22}&\ell_{23}
\end{pmatrix}
$$
where $\ell_{i,j}\in \mathbb{C}[x_0,\dotsc,x_4]_1$ with of course the 
condition $\lambda x_3+\mu x_4=0$. 
From this, eliminating the parameters $(\lambda:\mu)$, we can obtain a matrix 
$$
A_\eta:=\begin{pmatrix}
\ell_{11}&\ell_{12}&\ell_{13}\\
\ell_{21}&\ell_{22}&\ell_{23}\\
\ell_{31}&\ell_{32}&\ell_{33}
\end{pmatrix}
$$
whose determinant defines a cubic hypersurface which contains $F$---and 
in fact $C=V(\det(A_\eta))\cap \eta$. Therefore, $F$ is the degeneracy locus 
of the following $3\times 4$ matrix:
$$
A:=\begin{pmatrix}
\ell_{11}&\ell_{12}&\ell_{13}&x_3\\
\ell_{21}&\ell_{22}&\ell_{23}&x_4\\
\ell_{31}&\ell_{32}&\ell_{33}&0
\end{pmatrix}.
$$
The fact that $F$ generate a first order congruence follows from 
Theorem~2.6 of \cite{DP3} (in fact the argument of the proof is still valid 
even if $F$ is not smooth). 
Finally we conclude by Lemma~\ref{lem:bor}. 

\end{proof}

\begin{remark}
It is a well known fact that if the Bordiga surface 
$F$ is smooth, it is a blow-up
of $\p^2$ in $10$ points $P_1,\dotsc,P_{10}$ 
embedded in $\p^4$ by the linear system
$$
\abs{D}:=\abs{\pi^*4L-E_1-E_2-\dotsb-E_{10}}
$$  
where $\pi:F\rightarrow \p^2$ is the blow-up
in the points 
$P_1,\dotsc,P_{10}$, 
$L$ is a line in $\p^2$ and $E_i$ is the fibre of $\pi$ over $P_i$. 
In fact $F$ contains $10$
distinct lines and $10$ distinct plane cubics such that each line 
meets a single cubic. 

Even if $F$ is not smooth, one can show that $F$ is the blow-up of 
$\p^2$; with some work, one can give results similar to those for the 
Del Pezzo surfaces of \cite{dem}. For the sake of concision, we will not 
perform these calculations. 
\end{remark}

An example (due to F. Zak) of a singular Bordiga surface is the following 

\begin{example}
Consider the rational normal curve $C$ in $\p^5$, which is given by the 
$2\times 2$ minors of the catalecticant matrix
$$
M:=\begin{pmatrix}
x_{0}&x_{1}&x_{2}&x_{3}\\
x_{1}&x_{2}&x_{3}&x_{4}\\
x_{2}&x_{3}&x_{4}&x_5
\end{pmatrix};
$$
its secant variety $C^2$ is given by the $3\times 3$ minors of $M$: see 
\cite{hsv}. 
Then consider a general linear section $\ell$ of $C^2$: 
it is a Bordiga surface 
with $5$ singular points, which are $\ell\cap C$.   
\end{example}

From the results just proven in this section, we get 
an interesting characterization of 
the Bordiga surface:

\begin{proposition}
The only linearly normal surface 
in $\p^4$ whose trisecants generate a first 
order congruence is the Bordiga surface.
\end{proposition}

\begin{remark}
The Bordiga surface
is given by the vanishing of the 
minors of a general 
matrix of type $3\times 4$ of linear entries. 
We recall that also the rational 
normal cubic is the only curve in $\p^3$ whose secants generate a first 
order congruence, see \cite{DP1} and that this curve 
is given by the vanishing of the minors of a general matrix of 
type $2\times 3$ of linear entries. 
\end{remark}

\subsection{Final remarks on these congruences}
Let us make some concluding remarks about the results just proved: 
we consider now the congruence $B$ as a subvariety of dimension three in $\Gr(1,4)$. 
Let us give the following

\begin{proof}[Proof of Theorem~\textup{\ref{thm:smo}}]
First of all, the fact that congruence given by the trisecant lines of  
the Veronese surface is smooth follows from the general results of 
Subsection~2.1 of \cite{DP3}; 
similarly, the variety $B$ parametrising the trisecant 
lines of the Bordiga surface is smooth from the general results of 
Subsection~2.2 of \cite{DP3}, 
in particular Proposition~2.5 
and Theorem~2.6. 
For calculating the invariants of these congruences, we can simply get them 
from the resolutions 
of their ideal sheaves given in \cite{DP3}. 
Actually, the sectional genus $\pi(B)$ can be calculated also directly, and 
this is done, for example, in \cite{ArT}; in this 
case the formula is 
\begin{equation*}
\pi(B)=p_a(C_a)+a-1,
\end{equation*}
where $C_a$ is the curve in $\Gr(1,4)$ which corresponds to the scroll of the 
lines of the congruence contained in a hyperplane of $\p^4$. 
The genus $p_a(C_a)$ is calculated in Proposition~\textup{2.4} of \cite{GP} 
in terms of the degree and of the genus of $C_{H'}$, so we can compute the genus  
simply applying this result. But here we will compute $p_a(C_a)$ directly with a 
nice geometric argument. 
It is immediate that in the case of a linear congruence $C_a$ is rational; 
on the other hand, in the case of the Bordiga surface, 
we can consider a hyperplane section $H$ 
containing two twisted cubics $C_1$ and $C_2$ meeting in four points, 
$P_1,\dotsc,P_4$. 
The secants of $C_1$ which meet $C_2$ form, in the Grassmannian
$\Gr(1,H)\cong\Gr(1,3)$, a curve of degree $12$, 
as an easy computation with the Schubert
calculus shows. In this case, the two cubics meet in four points, so the 
curve of degree $12$ has, as irreducible components, four conics---given by the 
lines of the joins $P_i C_1$---which have to be cut out, and a quartic, $C_4$, 
which is rational, since it is contained in the Veronese surface (\ie the secant  
variety of $C_1$ in $\Gr(1,H)$).  
Using the same argument for $C_2$, we obtain another quartic $C'_4$, and so 
$C_a\cong C_4\cup C'_4$. $C_4\cap C'_4$ is given by four points: 
in fact, the Schubert calculus shows that there are $10$ lines which are secants both 
$C_1$ and $C_2$; but six of them are in the cones $P_i\cdot C_1$. 
Then, by adjunction, we have $p_a(C_a)=3$ and so $\pi(B)=10$.
\end{proof}

\begin{remark}
Note that the congruence $B$ given by the trisecant 
lines of a general projection 
to $\p^4$ of a rational normal scroll in $\p^6$ is not smooth: this is due 
to the fact that we have a $4$-parasitic plane, see Theorem~\ref{thm:m5}, and 
each of the lines of this plane is a singular (a $4$-tuple) point for $B$.
\end{remark}



\end{document}